\documentclass{IEEEtran4PSCC}
\ifCLASSINFOpdf
\else
\fi
\usepackage[cmex10]{amsmath}

\usepackage{graphicx}
\usepackage{fixltx2e}
\usepackage{stfloats}
\usepackage{tabularx}

\usepackage{epsfig} % for postscript graphics files
\usepackage{epstopdf}
\usepackage{amssymb}  % assumes amsmath package installed

\usepackage{cite}
\usepackage{url}
\usepackage{booktabs}
\usepackage[american,fulldiode]{circuitikz} %%Circuit drawing tools

\usepackage{algorithm}
\usepackage[noend]{algpseudocode}
\usepackage{float}
\newfloat{algorithm}{t}{lop}

\usepackage{tablefootnote}

\newcommand{\msosr}[1]{\mbox{\ensuremath{\text{MSOS}_{#1}}-$\mathbb{R}$}}
\newcommand{\msosc}[1]{\mbox{$\text{MSOS}_{#1}$-$\mathbb{C}$}}
\newcommand{\matpower}{\mbox{M{\sc atpower}}}

\allowdisplaybreaks

\begin{document}
\title{Computational Analysis of Sparsity-Exploiting Moment Relaxations of the OPF Problem}

%% To specify the authors when (number of affiliations <= 2)
\author{
\IEEEauthorblockN{Daniel K. Molzahn, Ian A. Hiskens}
\IEEEauthorblockA{%Electrical Engineering and Computer Science \\
EECS Department, University of Michigan\\
Ann Arbor, MI USA\\
\{molzahn,hiskens\}@umich.edu}
\and
\IEEEauthorblockN{C\'edric Josz, Patrick Panciatici}
\IEEEauthorblockA{R\&D Department, RTE \\
Paris, France\\
\{cedric.josz, patrick.panciatici\}@rte-france.com}
}

\maketitle

\begin{abstract}
With the potential to find global solutions, significant research
interest has focused on convex relaxations of the non-convex OPF
problem. Recently, ``moment-based'' relaxations from the Lasserre
hierarchy for polynomial optimization have been shown capable of
globally solving a broad class of OPF problems. Global solution of
many large-scale test cases is accomplished by exploiting sparsity and
selectively applying the computationally intensive higher-order
relaxation constraints. Previous work describes an iterative algorithm
that indicates the buses for which the higher-order constraints should
be enforced. In order to speed computation of the moment relaxations,
this paper provides a study of the key parameter in this algorithm as
applied to relaxations from both the original Lasserre hierarchy and a
recent complex extension of the Lasserre hierarchy.
\end{abstract}


% Use this to place sponsorships
\thanksto{The support of the Dow Sustainability Fellowship Program,
  ARPA-E grant \mbox{DE-AR0000232}, and Los Alamos National Laboratory
  subcontract 270958 is gratefully acknowledged.}

\section{Introduction}
The non-convex optimal power flow (OPF) problem is one of the most
important power system optimization problems. The OPF problem
minimizes an objective function subject to physical network
constraints and engineering limits. There is an extensive literature
of solution techniques for OPF problems, including successive
quadratic programs, Lagrangian relaxation, heuristic optimization, and
interior point methods~\cite{opf_litreview1993IandII,ferc4}. However,
while these local solution techniques often find global
solutions~\cite{molzahn_lesieutre_demarco-global_optimality_condition},
they may fail to converge or converge to a local
optimum~\cite{bukhsh_tps,ferc5}.

Recent research has focused on convex relaxations of the OPF problem. Convex relaxations lower bound the objective value and can certify infeasibility of OPF problems. For many OPF problems, a convex relaxation based on semidefinite programming (SDP) is \emph{exact} (i.e., the lower bound is tight and the solution provides the globally optimal decision variables)~\cite{lavaei_tps}. A second-order cone programming (SOCP) relaxation is provably exact for radial systems that satisfy certain non-trivial technical conditions~\cite{low_tutorial}. Developing tighter and faster relaxations is an active research area~\cite{sun2015,bienstock2015,coffrin2015}.

This paper focuses on a generalization of the SDP relaxation of~\cite{lavaei_tps} using the Lasserre hierarchy of semidefinite relaxations for real polynomial optimization problems~\cite{lasserre_book}. The primal form of the Lasserre hierarchy is interpreted as a truncated \emph{moment} series, while the dual form is interpreted as a \emph{sum-of-squares} optimization problem. We therefore refer to the \mbox{order-$\gamma$} relaxation in this hierarchy as \msosr{\gamma}. 

The OPF problem is a polynomial optimization problem in terms of the complex voltage phasors. Separating the complex voltages into real and imaginary parts yields a polynomial optimization problem in real variables. The first-order relaxation in the Lasserre hierarchy (i.e., \msosr{1}) is equivalent to the SDP relaxation of~\cite{lavaei_tps}. Higher-order moment relaxations globally solve a broader class of OPF problems at the computational cost of larger SDPs.

Recent work builds a complex hierarchy \msosc{\gamma} that is inspired
by the Lasserre hierarchy~\cite{complex_hierarchy}. Rather than
decomposing into real and imaginary parts, this complex hierarchy is
directly constructed from the complex voltages, which provides
computational advantages for many problems. This paper
summarizes \msosc{\gamma} through analogy to \msosr{\gamma}.

The computational requirements of both \msosr{\gamma} and \msosc{\gamma} grow quickly with increasing relaxation order $\gamma$. Fortunately, low-order relaxations ($\gamma\leqslant 2$) solve many small OPF problems~\cite{pscc2014,cedric,ibm_paper}. However, dense formulations for low-order relaxations are intractable for problems with more than approximately ten buses. Applying a ``chordal sparsity'' technique for the first-order relaxations enables solution of OPF problems with thousands of buses~\cite{jabr2011,molzahn_holzer_lesieutre_demarco-large_scale_sdp_opf}. Extending the chordal sparsity technique to higher-order relaxations~\cite{waki2006} facilitates the solution of the second-order relaxation for OPF problems with up to approximately 40 buses~\cite{molzahn_hiskens-sparse_moment_opf}.

Selectively enforcing the computationally intensive higher-order
constraints only at ``problematic'' buses enables extension to larger
problems~\cite{molzahn_hiskens-sparse_moment_opf}, including some OPF
problems with thousands of buses~\cite{cdc2015,complex_hierarchy}. The
iterative algorithm in~\cite{molzahn_hiskens-sparse_moment_opf}
applies the higher-order relaxation constraints to specific buses
identified with a heuristic, based on ``power injection mismatches'',
that provides an indication of the quality of the solution at each
bus. (Further details are provided in Section~\ref{l:sparsity}.)

The algorithm in~\cite{molzahn_hiskens-sparse_moment_opf} requires a
single parameter specifying how many buses should have higher-order
constraints applied at each iteration. There is a computational
trade-off associated with this parameter: applying higher-order
constraints at more buses potentially requires fewer iterations but a
larger computational cost per iteration. By balancing this trade-off,
this parameter can have a large impact on overall solution times.

To improve the computational effectiveness of both \msosr{\gamma} and \msosc{\gamma}, this paper presents a study of this parameter. Specifically, appropriate choice of this parameter is investigated by applying the algorithm in~\cite{molzahn_hiskens-sparse_moment_opf} to a variety of large test cases (Polish models in \matpower{}~\cite{matpower} and European models from the PEGASE project~\cite{pegase}).

This paper is organized as follows. Section~\ref{l:opf_formulation} describes the OPF problem in complex variables. Section~\ref{l:hierarchies} summarizes \msosr{\gamma} and \msosc{\gamma} as applied to the OPF problem, including the sparsity-exploiting algorithm in~\cite{molzahn_hiskens-sparse_moment_opf}. Section~\ref{l:computational_study} presents a computational study of the parameter in this algorithm. Section~\ref{l:conclusion} concludes the paper.

\section{OPF Problem Formulation}
\label{l:opf_formulation}

We first present an OPF formulation in terms of complex voltages, active and reactive power injections, and apparent-power line-flow limits. Consider an $n$-bus power system, where $\mathcal{N} = \left\lbrace 1, 2, \ldots, n \right\rbrace$ is the set of all buses, $\mathcal{G}$ is the set of generator buses, and $\mathcal{L}$ is the set of all lines. Let $P_{Dk} + \mathbf{j} Q_{Dk}$ represent the active and reactive load demand, where $\mathbf{j}$ is the imaginary unit, and $V_k$ the complex voltage phasor at each bus~$k \in \mathcal{N}$. Superscripts ``max'' and ``min'' denote specified upper and lower limits. Buses without generators have maximum and minimum generation set to zero. Let $\mathbf{Y}$ denote the network admittance matrix. The generator at bus~$k\in\mathcal{G}$ has a quadratic cost function for active power generation with coefficients $c_{k2} \geqslant 0$, $c_{k1}$, and $c_{k0}$. 

We use a line model with an ideal transformer that has a specified turns ratio $\tau_{lm} e^{\mathbf{j}\theta_{lm}} \colon 1$ in series with a $\Pi$ circuit with series impedance $R_{lm} + \mathbf{j} X_{lm}$ (equivalent to an admittance of $y_{lm} = \frac{1}{R_{lm} + \mathbf{j} X_{lm}}$) and shunt admittance $\mathbf{j} b_{sh,lm}$~\cite{matpower}.

The OPF problem is
\begin{subequations}
\label{eq:opf}
\begin{align}
\label{eq:opf_obj}
& \min_{V\in \mathbb{C}^n} \sum_{k \in \mathcal{G}} c_{k2}\left(V^H \mathbf{H}_k V \!+\! P_{Dk}\right)^2 \!+\! c_{k1}\left(V^H \mathbf{H}_k V \!+\! P_{Dk}\right) \!+\! c_{k0} \\[-9pt]
\nonumber & \text{subject to}
\end{align}
\vspace{-20pt}
\begin{align}
\label{eq:opfP}
& P_k^{\min} \leqslant V^H \mathbf{H}_k V + P_{Dk} \leqslant P_k^{\max} & \forall k \in \mathcal{N} \\
\label{eq:opfQ}
& Q_k^{\min} \leqslant V^H \tilde{\mathbf{H}}_k V + Q_{Dk} \leqslant Q_k^{\max} & \forall k \in \mathcal{N} \\
\label{eq:opfV}
& \left(V_k^{\min}\right)^2 \leqslant V^He_k^{\vphantom{\intercal}} e_k^\intercal V \leqslant \left(V_k^{\max}\right)^2 & \forall k \in \mathcal{N} \\
\nonumber
& \left(V^H\frac{\mathbf{F}_{lm}^{\vphantom{H}} \!+\! \mathbf{F}_{lm}^H}{2}V \right)^2 \!+\! \left(V^H\frac{\mathbf{F}_{lm}^{\vphantom{H}} \!-\! \mathbf{F}_{lm}^H}{2\mathbf{j}}V \right)^2 \!\leqslant \left(S_{lm}^{\max}\right)^2 \hspace{-150pt} \\ \label{eq:opfSlm} & & \forall \left(l,m\right) \in \mathcal{L} \\[-5pt]
\label{eq:opfVref}
&  V_1 - \overline{V}_1 = 0
\end{align}
\end{subequations}

% \nonumber
% & c_{k2}\left(V^H \mathbf{H}_k V \!+\! P_{Dk}\right)^2 \!+\! c_{k1}\left(V^H \mathbf{H}_k V \!+\! P_{Dk}\right) \!+\! c_{k0} \leqslant t_k\hspace{-100pt}\\\label{eq:opf_obj_constraint} & & \forall k \in \mathcal{G} \\

\noindent where $\left(\overline{\;\cdot\;}\right)$ denotes the complex conjugate, $\left(\cdot \right)^\intercal$ indicates the transpose, $\left(\cdot\right)^H$ indicates the complex conjugate transpose, and $e_k$ is the $k^{th}$ column of the identity matrix. Hermitian matrices in the active and reactive power injection constraints~\eqref{eq:opfP} and~\eqref{eq:opfQ}, respectively, are $\mathbf{H}_k = \frac{\mathbf{Y}^H e_k^{\vphantom{\intercal}} e_k^\intercal + e_k^{\vphantom{\intercal}} e_k^\intercal \mathbf{Y}}{2}$ and $\tilde{\mathbf{H}}_k = \frac{\mathbf{Y}^H e_k^{\vphantom{\intercal}} e_k^\intercal - e_k^{\vphantom{\intercal}} e_k^\intercal \mathbf{Y}}{2\mathbf{j}}$. Constraint~\eqref{eq:opfV} limits the squared voltage magnitude at bus~$k$. For the apparent-power line-flow constraints~\eqref{eq:opfSlm}, the Hermitian matrix $\mathbf{F}_{lm} = \left(\frac{\overline{y_{lm}} - \mathbf{j} b_{sh,lm}}{\tau_{lm}^2} \right) e_l^{\vphantom{\intercal}} e_l^\intercal + \left(\frac{-\overline{y_{lm}}}{\tau_{lm} e^{\mathbf{j}\theta_{lm}}} \right) e_m^{\vphantom{\intercal}} e_l^\intercal$. Flow constraints~\eqref{eq:opfSlm} are enforced at both line terminals ($\left(l,m\right) \in \mathcal{L}$ and  $\left(m,l\right) \in \mathcal{L}$). The angle reference is set by~\eqref{eq:opfVref}. Since $\mathbf{H}_k$, $\tilde{\mathbf{H}}_k$, $e_k^{\vphantom{\intercal}} e_k^\intercal$, and $\mathbf{F}_{lm}$ are Hermitian, all constraints in~\eqref{eq:opf} are real-valued polynomials in complex variables.

% , and~\eqref{eq:opf_obj_constraint} with auxiliary variable $t_k$ implements the quadratic cost function for the generator 

\section{Relaxation Hierarchies}
\label{l:hierarchies}

The OPF formulation is composed of polynomials, which enables global solution using polynomial optimization theory. Separating the complex decision variables into real and imaginary components $V_k = V_{dk} + \mathbf{j} V_{qk}$ facilitates the application of moment/sum-of-squares relaxations \msosr{\gamma} from the Lasserre hierarchy~\cite{lasserre_book} to the OPF problem~\cite{pscc2014,cedric,ibm_paper,molzahn_hiskens-sparse_moment_opf}. Directly building a relaxation hierarchy from the complex polynomial optimization formulation yields complex moment/sum-of-squares relaxations \msosc{\gamma}~\cite{complex_hierarchy}. This section reviews these hierarchies in the context of the OPF problem and discusses methods for exploiting sparsity. 

\subsection{Hierarchy for Real Polynomial Optimization Problems}
\label{l:real_hierarchy}

\begin{figure*}[!b]
\small
\setcounter{equation}{6}
\hrule
\begin{subequations}
\begin{align}\label{eq:realx2} \small
x_2 = & \left[\begin{array}{cccccccccc} 1 & V_{d1} & V_{d2} & V_{q2} & V_{d1}^2 & V_{d1}V_{d2} & V_{d1}V_{q2} & V_{d2}^2 & V_{d2}V_{q2} & V_{q2}^2\end{array}\right]^\intercal \\ \label{eq:compz2}
z_2 = & \left[\begin{array}{cccccc} 1 & V_{1} & V_{2} & V_{1}^2 & V_{1}V_{2} &  V_{2}^2\end{array}\right]^\intercal
\end{align}
\end{subequations}
\begin{subequations}
\noindent\begin{tabularx}{\textwidth}{@{}XX@{}}
{\vspace{-20pt}
\begin{align} \small \label{eq:realMoment2} 
& \mathbf{M}_2 \left\{y \right\} = L_y\left\{x_2^{\vphantom{\intercal}} x_2^\intercal\right\} = \\ \nonumber & \left[\begin{array}{c|ccc|cccccc} 
y_{000} & y_{100} & y_{010} & y_{001} & y_{200} & y_{110} & y_{101} & y_{020} & y_{011} & y_{002} \\\hline
y_{100} & y_{200} & y_{110} & y_{101} & y_{300} & y_{210} & y_{201} & y_{120} & y_{111} & y_{102} \\
y_{010} & y_{110} & y_{020} & y_{011} & y_{210} & y_{120} & y_{111} & y_{030} & y_{021} & y_{012} \\
y_{001} & y_{101} & y_{011} & y_{002} & y_{201} & y_{111} & y_{102} & y_{021} & y_{012} & y_{003} \\ \hline
y_{200} & y_{300} & y_{210} & y_{201} & y_{400} & y_{310} & y_{301} & y_{220} & y_{211} & y_{202} \\ 
y_{110} & y_{210} & y_{120} & y_{111} & y_{310} & y_{220} & y_{211} & y_{130} & y_{121} & y_{112} \\
y_{101} & y_{201} & y_{111} & y_{102} & y_{301} & y_{211} & y_{202} & y_{121} & y_{112} & y_{103} \\
y_{020} & y_{120} & y_{030} & y_{021} & y_{220} & y_{130} & y_{121} & y_{040} & y_{031} & y_{022} \\
y_{011} & y_{111} & y_{021} & y_{012} & y_{211} & y_{121} & y_{112} & y_{031} & y_{022} & y_{013} \\ 
y_{002} & y_{102} & y_{012} & y_{003} & y_{202} & y_{112} & y_{103} & y_{022} & y_{013} & y_{004} 
\end{array}\right]
\end{align}} & {\begin{align} \small \label{eq:compMoment2}
& \hat{\mathbf{M}}_2 \left\{\hat{y} \right\} = \hat{L}_{\hat{y}}\big\{z_2^{\vphantom{H}} z_2^H\big\} = \\ \nonumber & \left[\begin{array}{c|cc|ccc} 
\hat{y}_{00,00} & \hat{y}_{00,10} & \hat{y}_{00,01} & \hat{y}_{00,20} & \hat{y}_{00,11} & \hat{y}_{00,02} \\\hline
\hat{y}_{10,00} & \hat{y}_{10,10} & \hat{y}_{10,01} & \hat{y}_{10,20} & \hat{y}_{10,11} & \hat{y}_{10,02} \\
\hat{y}_{01,00} & \hat{y}_{01,10} & \hat{y}_{01,01} & \hat{y}_{01,20} & \hat{y}_{01,11} & \hat{y}_{01,02} \\\hline
\hat{y}_{20,00} & \hat{y}_{20,10} & \hat{y}_{20,01} & \hat{y}_{20,20} & \hat{y}_{20,11} & \hat{y}_{20,02} \\
\hat{y}_{11,00} & \hat{y}_{11,10} & \hat{y}_{11,01} & \hat{y}_{11,20} & \hat{y}_{11,11} & \hat{y}_{11,02} \\
\hat{y}_{02,00} & \hat{y}_{02,10} & \hat{y}_{02,01} & \hat{y}_{02,20} & \hat{y}_{02,11} & \hat{y}_{02,02}
\end{array}\right]
\end{align}}
\end{tabularx}
\end{subequations}
\vspace{0pt}
\begin{subequations}
\begin{align}\label{eq:realLocal2} \small \mathbf{M}_{1}\left\lbrace\left(f_{V2} - 0.81\right) y \right\rbrace & = \left[\begin{array}{c|ccc}
y_{020} + y_{002} - 0.81y_{000} & y_{120} + y_{102} - 0.81y_{100} & y_{030} + y_{012} - 0.81y_{010} & y_{021} + y_{003} - 0.81y_{001}  \\\hline
y_{120} + y_{102} - 0.81y_{100} & y_{220} + y_{202} - 0.81y_{200} & y_{130} + y_{112} - 0.81y_{110} & y_{121} + y_{103} - 0.81y_{101} \\
y_{030} + y_{012} - 0.81y_{010} & y_{130} + y_{112} - 0.81y_{110} & y_{040} + y_{022} - 0.81y_{020} & y_{031} + y_{013} - 0.81y_{011} \\
y_{021} + y_{003} - 0.81y_{001} & y_{121} + y_{103} - 0.81y_{101} & y_{031} + y_{013} - 0.81y_{011} & y_{022} + y_{004} - 0.81y_{002}  \\
\end{array}\right] \\[-8pt]
\label{eq:compLocal2} 
\hat{\mathbf{M}}_{1}\left\lbrace\left(\hat{f}_{V2} - 0.81\right) \hat{y} \right\rbrace & = \left[\begin{array}{c|ccc}
\hat{y}_{01,01}-0.81\hat{y}_{00,00} & \hat{y}_{01,11}-0.81\hat{y}_{00,10} & \hat{y}_{01,02}-0.81\hat{y}_{00,01} \\ \hline
\hat{y}_{11,01}-0.81\hat{y}_{10,00} & \hat{y}_{11,11}-0.81\hat{y}_{10,10} & \hat{y}_{11,02}-0.81\hat{y}_{10,01} \\
\hat{y}_{02,01}-0.81\hat{y}_{01,00} & \hat{y}_{02,11}-0.81\hat{y}_{01,10} & \hat{y}_{02,02}-0.81\hat{y}_{01,01}
\end{array}\right]
\end{align}
\end{subequations}
\setcounter{equation}{1}
\end{figure*}

The Lasserre hierarchy \msosr{\gamma} builds relaxations that take the form of semidefinite programs. For generic polynomial optimization problems that satisfy a technical condition\footnote{This technical condition is satisfied when the decision variables are bounded and is therefore not restrictive for the OPF problem.}, relaxations in the Lasserre hierarchy converge to the global solutions of generic polynomial optimization problems at a finite relaxation order~\cite{lasserre_book,nie2014}. While the Lasserre hierarchy can find all global solutions, we focus on OPF problems with a single global optimum.\footnote{See~\cite{molzahn_hiskens-sparse_moment_opf} for a discussion of atypical cases with multiple global optima.}

We begin with several definitions. Define the vector of real decision
variables $\xi \in \mathbb{R}^{2n}$ as $\xi := \begin{bmatrix} V_{d1}
  & V_{d2} & \ldots V_{qn} \end{bmatrix}^\intercal$.\footnote{The
  ability to select an angle reference in the OPF problem enables
  specification of one arbitrarily chosen variable. We choose $V_{q1}
  = 0$.} A monomial is defined using a vector $\alpha \in
\mathbb{N}^{2n}$ of exponents: $\xi^\alpha := V_{d1}^{\alpha_1}
V_{d2}^{\alpha_2}\cdots V_{qn}^{\alpha_{2n}}$. A polynomial is
$g\left(\xi\right) := \sum_{\alpha \in \mathbb{N}^{2n}} g_{\alpha}
\xi^{\alpha}$, where $g_{\alpha}$ is the real scalar coefficient
corresponding to the monomial $\xi^{\alpha}$.

Define a linear functional $L_y\left\lbrace g\right\rbrace$ which
replaces the monomials $\xi^{\alpha}$ in a polynomial
$g\left(\xi\right)$ with real scalar variables $y$:
\begin{equation}
\label{eq:Lreal}
L_y\left\lbrace g \right\rbrace := \sum_{\alpha \in \mathbb{N}^{2n}} g_{\alpha} y_{\alpha}.
\end{equation}
For a matrix $g\left(\xi\right)$, $L_y\left\lbrace g\right\rbrace$ is applied componentwise. 

Consider, for example, the vector $\xi = \begin{bmatrix}V_{d1} & V_{d2} & V_{q2} \end{bmatrix}^\intercal$ corresponding to the voltage components of a two-bus system, where the angle reference is used to eliminate $V_{q1}$, and the polynomial $g\left(\xi\right) = -\left(0.9\right)^2 + V_{d2}^2 + V_{q2}^2$. (The constraint $g\left(\xi\right) \geqslant 0$ forces the voltage magnitude at bus~2 to be greater than or equal to 0.9~per unit.) Then $L_y\left\lbrace g\right\rbrace = -\left(0.9\right)^2y_{000} + y_{020} + y_{002}$. Thus, $L_y\left\lbrace g \right\rbrace$ converts a polynomial $g\left(\xi\right)$ to a linear function of $y$.

For the order-$\gamma$ relaxation, define a vector $x_\gamma$
consisting of all monomials of the voltage components up to order
$\gamma$:
\begin{align}
\nonumber
x_\gamma := & \left[ \begin{array}{ccccccc} 1 & V_{d1} & \ldots & V_{qn} & V_{d1}^2 & V_{d1}V_{d2} & \ldots \end{array} \right. \\* \label{eq:x_d}
& \qquad \left.\begin{array}{cccccc} \ldots & V_{qn}^2 & V_{d1}^3 & V_{d1}^2 V_{d2} & \ldots & V_{qn}^\gamma \end{array}\right]^\intercal.
\end{align}

The relaxations are composed of positive semidefinite constraints on \emph{moment} and \emph{localizing} matrices. The symmetric moment matrix $\mathbf{M}_{\gamma}$ is composed of entries $y_\alpha$ corresponding to all monomials $\xi^{\alpha}$ up to order $2\gamma$:

\begin{equation}
\label{eq:real_moment}
\mathbf{M}_\gamma \left\lbrace y \right\rbrace := L_y\left\lbrace x_\gamma^{\vphantom{\intercal}} x_\gamma^\intercal\right\rbrace.
\end{equation}

Symmetric localizing matrices are defined for each constraint of~\eqref{eq:opf}. For a polynomial constraint $g\left(\xi\right) \geqslant 0$ of degree $2\eta$, the localizing matrix is:
\begin{equation}
\label{eq:real_local}
\mathbf{M}_{\gamma - \eta} \left\lbrace g y \right\rbrace := L_y \left\lbrace g x_{\gamma-\eta}^{\vphantom{\intercal}} x_{\gamma-\eta}^{\intercal} \right\rbrace.
\end{equation}
See~\eqref{eq:realMoment2} and~\eqref{eq:realLocal2} for example moment and localizing matrices for \msosr{2} applied to a two-bus OPF problem. 

Separating real and imaginary parts $V_k = V_{dk} + \mathbf{j} V_{qk}$, define $f_{Pk}$, $f_{Qk}$, and $f_{Vk}$ for the active power injection, reactive power injection, and squared voltage magnitudes at bus~$k$ using~\eqref{eq:opfP}, \eqref{eq:opfQ}, and \eqref{eq:opfV}, respectively, as functions of the real variables $V_d$ and $V_q$. Likewise, define $f_{Slm}$ for the squared apparent power flow on the line from bus~$l$ to bus~$m$ using~\eqref{eq:opfSlm} and $f_{Ck}$ as the quadratic cost function for the generator at bus~$k$ using \eqref{eq:opf_obj} as functions of the real variables $V_d$ and $V_q$.

%\begin{subequations}
%\begin{align}
%f_{Pk} & := \left(V_{d} + \mathbf{j} V_q\right)^H \mathbf{H}_k \left(V_{d} + \mathbf{j} V_q\right) + P_{Dk} \\
%f_{Qk} & := \left(V_{d} + \mathbf{j} V_q\right)^H \tilde{\mathbf{H}}_k \left(V_{d} + \mathbf{j} V_q\right) + Q_{Dk} \\
%f_{Vk} & := \left(V_{d} + \mathbf{j} V_q\right)^H e_k^{\vphantom{\intercal}} e_k^\intercal \left(V_{d} + \mathbf{j} V_q\right) \\
%f_{Slm} & := 
%\end{align}
%\end{subequations}

The order-$\gamma$ relaxation \msosr{\gamma} is:
\begin{subequations}
\label{eq:msosr}
\begin{align}
\label{eq:msosr_obj}& \min_{y} L_y\left\lbrace \sum_{k \in \mathcal{G}} f_{Ck} \right\rbrace \qquad \mathrm{subject\; to} \hspace{-150pt} &  \\
\label{eq:msosr_Pmin} &  \mathbf{M}_{\gamma-1}\left\lbrace
\left(f_{Pk} - P_k^{\min}\right) y \right\rbrace \succcurlyeq 0 &
\forall k\in\mathcal{N} \\
%\end{align}\begin{align} %\\
\label{eq:msosr_Pmax} &  \mathbf{M}_{\gamma-1}\left\lbrace \left(P_k^{\max} - f_{Pk} \vphantom{P_k^{\min}}\right) y \right\rbrace \succcurlyeq 0 & \forall k\in\mathcal{N}\\
\label{eq:msosr_Qmin} &  \mathbf{M}_{\gamma-1}\left\lbrace \left(f_{Qk} - Q_k^{\min}\right) y \right\rbrace \succcurlyeq 0 & \forall k\in\mathcal{N}\\
\label{eq:msosr_Qmax} &  \mathbf{M}_{\gamma-1}\left\lbrace \left(Q_k^{\max} - f_{Qk}  \vphantom{P_k^{\min}}\right) y \right\rbrace \succcurlyeq 0 & \forall k\in\mathcal{N}\\
\label{eq:msosr_Vmin} &  \mathbf{M}_{\gamma-1}\left\lbrace \left(f_{Vk} - \left(V_k^{\min}\right)^2\right) y \right\rbrace \succcurlyeq 0 & \forall k\in\mathcal{N}\\
\label{eq:msosr_Vmax} &  \mathbf{M}_{\gamma-1}\left\lbrace \left(\left(V_k^{\max}\right)^2 - f_{Vk}  \vphantom{P_k^{\min}}\right) y \right\rbrace \succcurlyeq 0 & \forall k\in\mathcal{N} \\
\label{eq:msosr_Smax} & \mathbf{M}_{\gamma-2}\left\lbrace \left( \left(S_{lm}^{\max}\right)^2 - f_{Slm} \vphantom{P_k^{\min}}\right) y \right\rbrace \succcurlyeq 0 & \forall \left(l,m\right)\in\mathcal{L} \\
\label{eq:msosr_Msdp} & \mathbf{M}_\gamma \left\{y\right\} \succcurlyeq 0 & \\
\label{eq:msosr_y0} & y_{00\ldots 0} = 1 & \\
\label{eq:msosr_Vref} & y_{\star\star\ldots\star\rho\star\ldots\star} = 0 & \rho = 1,\ldots,2\gamma.
\end{align}
\end{subequations}

\noindent where $\succcurlyeq 0$ indicates that the corresponding matrix is positive semidefinite and $\star$ represents any integer in $\left[ 0,\; 2\gamma -1\right]$. The constraint~\eqref{eq:msosr_y0} enforces the fact that $x^{0} = 1$. The constraint~\eqref{eq:msosr_Vref} corresponds to the angle reference $V_{q1} = 0$; the $\rho$ in \eqref{eq:msosr_Vref} is in the index $n+1$, which corresponds to the variable $V_{q1}$. Note that the angle reference can alternatively be used to eliminate all terms corresponding to $V_{q1}$ to reduce the size of the semidefinite program.

A dual form of the ``moment'' relaxation presented here is a \emph{sum-of-squares} program, thus leading to the nomenclature \msosr{\gamma}~\cite{lasserre_book}. There is zero duality gap between the moment and sum-of-squares formulations for the OPF problem~\cite{josz-2015}.

The order-$\gamma$ relaxation yields a single global solution if
$\mathrm{rank}\left(\mathbf{M}_{\gamma}\left\{y\right\}\right) = 1$. The
global solution $V^\ast$ to the OPF problem~\eqref{eq:opf} is then
determined by a spectral decomposition of the diagonal block of the
moment matrix corresponding to the second-order monomials (i.e.,
$\left|\alpha\right| = 2$, where $\left|\;\cdot\; \right|$ indicates
the one-norm). Specifically, let $\eta$ be a unit-length eigenvector
corresponding to the non-zero eigenvalue $\lambda$ of
$\left[\mathbf{M}_1 \{y\}\right]_{\left(2:k,2:k\right)}$, where $k =
2n+1$ and subscripts indicate the vector entries in MATLAB
notation.\footnote{For the example
  \eqref{eq:realx2}, \eqref{eq:realMoment2}, the angle reference was
  established by eliminating $V_{q1}$. Therefore, $k=2n$ for this
  case.}  Then the vector $V^\ast = \sqrt{\lambda} \left(\eta_{1:n} +
\mathbf{j} \eta_{\left(n+1\right):2n}\right)$ is the globally optimal
voltage phasor vector.

Note that the order $\gamma$ of the relaxation must be greater than or equal to half of the degree of any polynomial in the OPF problem~\eqref{eq:opf}. Although direct implementation of~\eqref{eq:opf} requires $\gamma \geqslant 2$ due to the fourth-order polynomials in the cost function~\eqref{eq:opf_obj} and apparent-power line-flow constraints~\eqref{eq:opfSlm}, these can be rewritten using a Schur complement~\cite{lavaei_tps} to allow $\gamma \geqslant 1$. Experience suggests that implementing~\eqref{eq:opf_obj} and~\eqref{eq:opfSlm} both directly and with a Schur complement formulation gives superior results for $\gamma \geqslant 2$.

\subsection{Hierarchy for Complex Polynomial Optimization Problems}
\label{l:complex_hierarchy}

Rather than first separating the decision variables into real and imaginary parts $V_d$ and $V_q$, a hierarchy of relaxations, denoted as \msosc{\gamma}, built directly from~\eqref{eq:opf} in complex variables has computational advantages for many OPF problems~\cite{complex_hierarchy}. Presentation of the complex hierarchy \msosc{\gamma} mirrors the development of \msosr{\gamma} in Section~\ref{l:real_hierarchy}.

We again begin with several definitions. Define the vector of complex decision variables $\zeta \in \mathbb{C}^n$ as $\zeta := \begin{bmatrix} V_1 & V_2 & \ldots & V_n \end{bmatrix}^\intercal$. A complex monomial is defined using two vectors of exponents $\alpha,\, \beta \in \mathbb{N}^n$: $\zeta^{\alpha} \overline{\zeta}^{\beta} := V_1^{\alpha_1}\cdots V_n^{\alpha_n} \overline{V_1}^{\beta_1} \cdots \overline{V_n}^{\beta_n}$. A polynomial $g\left(\zeta\right) := \sum_{\alpha,\beta \in \mathbb{N}^n} g_{\alpha,\beta} \zeta^{\alpha} \overline{\zeta}^{\beta}$, where $g_{\alpha,\beta}$ is the complex scalar coefficient corresponding to the monomial $\zeta^\alpha \overline{\zeta}^\beta$. Since $g\left(\zeta\right)$ outputs a real value, $\overline{g_{\alpha,\beta}} = g_{\beta,\alpha}$.

\setcounter{equation}{9}

Define a linear functional $\hat{L}_{\hat{y}}\left(g\right)$ which replaces the monomials in a polynomial $g\left(\zeta\right)$ with complex scalar variables $\hat{y}$:
\begin{equation}
\label{eq:Lcomp}
\hat{L}_{\hat{y}}\left\lbrace g \right\rbrace := \sum_{\alpha,\beta \in \mathbb{N}^{n}} g_{\alpha,\beta} \hat{y}_{\alpha,\beta}.
\end{equation}
For a matrix $g\left(\zeta\right)$, $\hat{L}_{\hat{y}}\left\lbrace g\right\rbrace$ is applied componentwise. 

Consider, for example, the vector $\zeta = \begin{bmatrix}V_{1} & V_{2} \end{bmatrix}^\intercal$ corresponding to the complex voltage phasors of a two-bus system and the polynomial $g\left(\zeta\right) = -\left(0.9\right)^2 + V_2 \overline{V_{2}}$. (The constraint $g\left(\zeta\right) \geqslant 0$ forces the voltage magnitude at bus~2 to be greater than or equal to 0.9~per unit.) Then $\hat{L}_{\hat{y}}\left\lbrace g\right\rbrace = -\left(0.9\right)^2\hat{y}_{00,00} + \hat{y}_{01,01}$. Thus, $\hat{L}_{\hat{y}}\left\lbrace g \right\rbrace$ converts a polynomial $g\left(\zeta\right)$ to a linear function of $\hat{y}$.

For the order-$\gamma$ relaxation, define a vector $z_\gamma$ consisting of all monomials of the voltages up to order $\gamma$ without complex conjugate terms (i.e., $\beta = 00\cdots 0$):
\begin{align}
\nonumber
z_\gamma := & \left[ \begin{array}{ccccccc} 1 & V_{1} & \ldots & V_{n} & V_{1}^2 & V_{1}V_{2} & \ldots \end{array} \right. \\ \label{eq:z_d}
& \qquad \left.\begin{array}{cccccc} \ldots & V_{n}^2 & V_{1}^3 & V_{1}^2 V_{2} & \ldots & V_{n}^\gamma \end{array}\right]^\intercal.
\end{align}

For the complex hierarchy, the symmetric moment matrix $\hat{\mathbf{M}}_{\gamma}$ is composed of entries $y_{\alpha,\beta}$ corresponding to all monomials $\zeta^{\alpha}\overline{\zeta}^{\beta}$ such that $\left|\alpha\right| + \left|\beta\right| \leqslant 2\gamma$:
\begin{equation}
\label{eq:comp_moment}
\hat{\mathbf{M}}_\gamma \left\lbrace \hat{y} \right\rbrace := \hat{L}_{\hat{y}}\left\lbrace z_\gamma^{\vphantom{H}} z_\gamma^H\right\rbrace.
\end{equation}

Symmetric localizing matrices are defined for each constraint
of~\eqref{eq:opf}. For a polynomial constraint $g\left(\zeta\right)
\geqslant 0$ with largest degree $\left|\alpha + \beta \right|$ among
all monomials equal to $2\eta$,\footnote{For OPF problems,
  all constraint and cost-function polynomials (other than the reference angle~\eqref{eq:opfVref}, which is discussed later) have even degree.} the
localizing matrix is:
\begin{equation}
\label{eq:comp_local}
\hat{\mathbf{M}}_{\gamma - \eta} \left\lbrace g \hat{y} \right\rbrace := \hat{L}_{\hat{y}} \left\lbrace g z_{\gamma-\eta}^{\vphantom{H}} z_{\gamma-\eta}^{H} \right\rbrace.
\end{equation}
See~\eqref{eq:compMoment2} and~\eqref{eq:compLocal2} for example moment and localizing matrices for \msosc{2} applied to a two-bus OPF problem. 

Define $\hat{f}_{Pk}$, $\hat{f}_{Qk}$, and $\hat{f}_{Vk}$ for the active power injection, reactive power injection, and squared voltage magnitudes at bus~$k$ using~\eqref{eq:opfP}, \eqref{eq:opfQ}, and \eqref{eq:opfV}, respectively, as functions of complex variables $V$. Likewise, define $\hat{f}_{Slm}$ for the squared apparent power flow on the line from bus~$l$ to bus~$m$ using~\eqref{eq:opfSlm} and $\hat{f}_{Ck}$ as the quadratic cost function for the generator at bus~$k$ using \eqref{eq:opf_obj} as functions of the complex variables $V$.

The order-$\gamma$ relaxation \msosc{\gamma} is
\begin{subequations}
\label{eq:msosc}
\begin{align}
\label{eq:msosc_obj} & \min_{\hat{y}} \hat{L}_{\hat{y}}\left\lbrace \sum_{k \in \mathcal{G}} \hat{f}_{Ck} \right\rbrace \qquad \mathrm{subject\; to} \hspace{-150pt} &  \\
\label{eq:msosc_Pmin} & \hat{\mathbf{M}}_{\gamma-1}\left\lbrace \left(\hat{f}_{Pk} - P_k^{\min}\right) \hat{y} \right\rbrace \succcurlyeq 0 & \forall k\in\mathcal{N}\\
\label{eq:msosc_Pmax} &  \hat{\mathbf{M}}_{\gamma-1}\left\lbrace \left(P_k^{\max} - \hat{f}_{Pk} \vphantom{P_k^{\min}}\right) \hat{y} \right\rbrace \succcurlyeq 0 & \forall k\in\mathcal{N}\\
\label{eq:msosc_Qmin} & \hat{\mathbf{M}}_{\gamma-1}\left\lbrace \left(\hat{f}_{Qk} - Q_k^{\min}\right) \hat{y} \right\rbrace \succcurlyeq 0 & \forall k\in\mathcal{N}\\
\label{eq:msosc_Qmax} &  \hat{\mathbf{M}}_{\gamma-1}\left\lbrace \left(Q_k^{\max} - \hat{f}_{Qk}  \vphantom{P_k^{\min}}\right) \hat{y} \right\rbrace \succcurlyeq 0 & \forall k\in\mathcal{N}\\
\label{eq:msosc_Vmin} &  \hat{\mathbf{M}}_{\gamma-1}\left\lbrace \left(\hat{f}_{Vk} - \left(V_k^{\min}\right)^2\right) \hat{y} \right\rbrace \succcurlyeq 0 & \forall k\in\mathcal{N}\\
\label{eq:msosc_Vmax} &  \hat{\mathbf{M}}_{\gamma-1}\left\lbrace \left(\left(V_k^{\max}\right)^2 - \hat{f}_{Vk}  \vphantom{P_k^{\min}}\right) \hat{y} \right\rbrace \succcurlyeq 0 & \forall k\in\mathcal{N} \\
\label{eq:msosc_Smax} & \hat{\mathbf{M}}_{\gamma-2}\left\lbrace \left(\left(S_{lm}^{\max}\right)^2 - \hat{f}_{Slm} \vphantom{P_k^{\min}}\right) \hat{y} \right\rbrace \succcurlyeq 0 & \forall \left(l,m\right)\in\mathcal{L} \\
\label{eq:msosc_Msdp} & \hat{\mathbf{M}}_\gamma \{\hat{y}\} \succcurlyeq 0 & \\
\label{eq:msosc_y0} & \hat{y}_{0\ldots 0,0\ldots 0} = 1.
\end{align}
\end{subequations}

\noindent Rather than explicitly setting an angle reference, we rotate the solution to~\eqref{eq:msosc} in order to satisfy~\eqref{eq:opfVref}. Note that a Schur complement formulation for~\eqref{eq:opf_obj} and~\eqref{eq:opfSlm} enables solution of~\eqref{eq:msosc} with $\gamma = 1$~\cite{lavaei_tps}. 

Similar to the real hierarchy, \msosc{\gamma} yields a single global
solution if
$\mathrm{rank}\big(\hat{\mathbf{M}}_{\gamma}\{\hat{y}\}\big)
= 1$. The global solution $V^\ast$ is calculated using a spectral
decomposition of the diagonal block of the moment matrix corresponding
to the second-order monomials (i.e., $\left|\alpha\right| =
\left|\beta\right| = 1$). Let $\hat{\eta}$ be a unit-length
eigenvector corresponding to the non-zero eigenvalue $\hat{\lambda}$
of $\big[\hat{\mathbf{M}}_1 \{\hat{y}\}
  \big]_{\left(2:n+1,2:n+1\right)}$. Then the vector $V^\ast =
\sqrt{\hat{\lambda}} \hat{\eta}$, rotated to match the angle
reference, gives the globally optimal voltages.

We next summarize theoretical developments for \msosc{\gamma} applied
to the OPF problem.\footnote{These statements do not hold for all
  polynomial optimization problems.} These developments, which are
proven in~\cite{complex_hierarchy}, relate to duality, convergence, and
comparison to the real hierarchy. See~\cite{lasserre_book} for similar
results for \msosr{\gamma}.

\begin{enumerate}
\item
Analogous to the real hierarchy, \msosc{\gamma} can be interpreted in a dual sum-of-squares form. There is zero duality gap between the primal~\eqref{eq:msosc} and  dual forms.
\item
The relaxations \msosr{1} and \msosc{1}, and the relaxation in~\cite{lavaei_tps} all give the same optimal objective values.
\item
Augmenting the OPF problem with a \emph{sphere constraint} 

\begin{equation}
\label{eq:sphere_constraint}
\sum_{i=1}^n \left(V_i \overline{V_i}\right) + \psi \overline{\psi} = \sum_{i=1}^n \left(V_i^{\max}\right)^2,
\end{equation}

\noindent where $\psi$ is a slack variable, guarantees convergence of \msosc{\gamma} to the global optimum with increasing $\gamma$. Observe that the sphere constraint is redundant due to the upper voltage magnitude limits~\eqref{eq:opfV} and therefore does not impair the applicability of \msosc{\gamma}.

\item
For \msosc{\gamma} with the sphere constraint~\eqref{eq:sphere_constraint} and \msosr{\gamma} of the same order $\gamma$:
\begin{enumerate}
\item
The optimal objective value from \msosr{\gamma} is at least as large as that from \msosc{\gamma}.

\item There exist optimization problems for which \msosr{\gamma} gives strictly superior objective values.

\item
Numerical results support the conjecture (formalized
in~\cite{complex_hierarchy}) that \msosc{\gamma} yields the same
optimal objective value as \msosr{\gamma} when applied to complex
polynomial optimization problems that exhibit rotational
  symmetry, such as OPF problems.
\end{enumerate}

\item
For an $n$-bus system, the size of the moment matrices~\eqref{eq:msosr_Msdp} and~\eqref{eq:msosc_Msdp} for \msosr{\gamma} (using the angle reference~\eqref{eq:msosr_Vref} to eliminate $V_{q1}$) and \msosc{\gamma} (converted to real representation for input to the solver~\cite[Example 4.42]{boyd2009}) are $\left(2n-1+\gamma\right)! / \left( \left(2n-1\right)! \gamma!\right)$ and $2\left(\left(n+\gamma\right)!\right)/\left(n!\gamma!\right)$, respectively. For example, $n=10$ and $\gamma=3$ correspond to matrices of size $1,\!540\times 1,\!540$ and $572\times 572$ for the real and complex hierarchies, respectively. Thus, the complex hierarchy is significantly more computationally tractable than the real hierarchy.
% n = 10; gamma = 3; msosr = factorial(2*n-1+gamma)/(factorial(2*n-1)*factorial(gamma)); msosc = 2*factorial(n+gamma)/(factorial(n)*factorial(gamma)); disp([msosr msosc]);
\end{enumerate}

\subsection{Exploiting Network Sparsity}
\label{l:sparsity}

Although \msosc{\gamma} is computationally superior to \msosr{\gamma}, matrices for both hierarchies grow quickly with the relaxation order $\gamma$. Practical evaluation of both hierarchies for large problems requires low relaxation order. Fortunately, low-order relaxations often yield high-quality lower bounds and, in many cases, the global optima for practical OPF problems.

However, for $n \gtrsim 10$, the ``dense'' formulations~\eqref{eq:msosr} and~\eqref{eq:msosc} are intractable even for $\gamma = 2$. Exploiting network sparsity enables the application of the hierarchies to many large problems. First proposed for the OPF problem in~\cite{jabr2011}, a ``chordal sparsity'' technique enables solution of the first-order relaxation for systems with thousands of buses~\cite{molzahn_holzer_lesieutre_demarco-large_scale_sdp_opf}. The approach in~\cite{waki2006} extends these chordal sparsity techniques to higher relaxation orders of the real hierarchy and can be readily applied to the complex hierarchy. Each positive semidefinite matrix constraint is replaced by a set of positive semidefinite matrix constraints on certain submatrices. These submatrices are defined by the maximal cliques of a specifically constructed chordal extension of the network graph.  See~\cite{molzahn_hiskens-sparse_moment_opf,complex_hierarchy} for a detailed description of this approach. 

When the relaxation order $\gamma=1$, the ``sparse'' versions of the moment/sum-of-squares hierarchies give equivalent solutions to their dense counterparts. However, this is not the case for $\gamma > 1$: the sparse hierarchies are generally not as tight as the dense hierarchies~\cite{lasserre_book}. Nevertheless, low-order sparse hierarchies globally solve many practical problems.

Applying the chordal sparsity approach enables solution of OPF problems with $n \lesssim 40$. Extension to larger problems is possible by recognizing that the computationally challenging higher-order relaxations are only necessary for the constraints associated with specific buses. In other words, rather than a single relaxation order $\gamma$ applied to all buses in the problem, each bus~$i$ has an associated relaxation order $\gamma_i$. By both exploiting sparsity and \emph{selectively applying} the computationally intensive higher-order relaxation constraints, many large OPF problems are computationally tractable.

Selectively applying the higher-order constraints requires a method
for determining $\gamma_i$ for each bus. We use a heuristic based on
``power injection mismatches'' to the closest rank-one
matrix~\cite{molzahn_hiskens-sparse_moment_opf}. Consider the
application of the complex hierarchy to an OPF problem with a single
global optimum. (Application of the real hierarchy proceeds
analogously.) An ``approximate'' solution $z^{\text{approx}}$ to the
OPF problem can be obtained from the largest eigenvalue
$\hat{\lambda}_1$ and associated unit-length eigenvector
$\hat{\eta}_1$ of the matrix
$\big[\hat{\mathbf{M}}_1\{\hat{y}\}\big]_{\left(2:n+1,2:n+1\right)}$,
so $z^{\text{approx}} := \sqrt{\hat{\lambda}_1}\,\hat{\eta}_1$. For
each bus~$i$, define a power injection mismatch $S_i^{\text{mis}}$
between the solution to the relaxation and $z^{\text{approx}}$:
\begin{align}\nonumber
S_i^{\text{mis}} := & \left| \left(\hat{f}_{Pi}\left(z^{\text{approx}}\right) - \hat{L}_{\hat{y}}\left\lbrace f_{Pi}\right\rbrace\right) \right. \\\label{eq:mismatch}
 & \qquad \left. + \; \mathbf{j} \left(\hat{f}_{Qi}\left(z^{\text{approx}}\right) - \hat{L}_{\hat{y}}\left\lbrace f_{Qi}\right\rbrace \right) \right|.
\end{align}

\begin{algorithm}
\caption{Iterative Solution for Sparse Relaxations}\label{alg:heuristic}
\begin{algorithmic}[1]
\State Initialize $\gamma_i := 1,\, i=1,\ldots, n$
\Repeat
\State Solve relaxation with order $\gamma$
\State Calculate mismatches $S_i^{\text{mis}},\, i=1,\ldots,n$ using~\eqref{eq:mismatch}
\State Increase entries of $\gamma$ according to mismatch heuristic
\Until{$\left|S_i^{\text{mis}}\right|_\infty < \epsilon $}
\State Extract solution $V^{\ast}$
\end{algorithmic}
\end{algorithm}

We employ an iterative algorithm for determining relaxation orders $\gamma_i,\, i=1,\ldots,n$. (See Algorithm~\ref{alg:heuristic}.) Each iteration of the algorithm solves the relaxation after increasing the relaxation orders $\gamma_i$ in a manner that is dependent on the largest associated $S_i^{\text{mis}}$ values. At each iteration of the algorithm, calculate $\gamma^{\text{max}} := \max_i \left\lbrace \gamma_i \right\rbrace$.\footnote{Note that $\gamma^{\text{max}}$ is not a specified maximum but can change at each iteration.} Each iteration increments $\gamma_i$ at up to $h$ buses, where $h$ is a specified parameter, that have the largest mismatches $S_i^{\text{mis}}$ among all buses satisfying two conditions: 1.)~$\gamma_i < \gamma^{\text{max}}$ and 2.)~$S_i^{\text{mis}} > \epsilon$, where $\epsilon$ is a specified mismatch tolerance. If no buses satisfy both of these conditions, increment $\gamma_i$ at up to $h$ buses with the largest $S_i^{\text{mis}}$ greater than the specified tolerance and increment $\gamma^{\text{max}}$. That is, in order to avoid unnecessarily increasing the size of the positive semidefinite matrices, the heuristic avoids incrementing the maximum relaxation order $\gamma^{\text{max}}$ until $\gamma_i = \gamma^{\text{max}}$ at all buses with mismatch $S_i^{\text{mis}} > \epsilon$. The algorithm terminates when $\left|S_i^{\text{mis}}\right|_{\infty} \leqslant \epsilon$, where $\left|\;\cdot\;\right|_{\infty}$ denotes the maximum absolute value, which indicates satisfaction of the rank condition for practical purposes. Thus, the relaxations are successively tightened in a manner that preserves computational tractability. 

%Since the Algorithm~\ref{alg:heuristic} eventually proceeds to build the complete moment/sum-of-squares hierarchies, Algorithm~\ref{alg:heuristic} inherits the theoretical convergence guarantees for \msosr{\gamma} and \msosc{\gamma}.

\section{Computational Study of the Parameter in the Sparsity-Exploiting Algorithm}
\label{l:computational_study}

There is a computational trade-off in choosing the value of $h$. Larger values of $h$ may result in fewer iterations of the algorithm but each iteration is slower if more buses than necessary have high-order relaxations. Smaller values of $h$ result in faster solution at each iteration, but may require more iterations. Previous work~\cite{molzahn_hiskens-sparse_moment_opf,cdc2015,complex_hierarchy} chose $h = 2$ based on limited computational experience. This section presents a more comprehensive study of the impact of this parameter. 

In theory, the choice of $h$ has no impact on the eventual convergence of the hierarchies. With sufficiently small tolerance $\epsilon$, Algorithm~\ref{alg:heuristic} eventually proceeds to build the complete hierarchies. Thus, Algorithm~\ref{alg:heuristic} inherits the theoretical convergence guarantees for \msosr{\gamma} and \msosc{\gamma}. In practice, the choice of $h$ significantly affects the computational requirements, and thus the practical capabilities of the hierarchies.

To study the impact of the parameter $h$ on solution times, we solved
both \msosr{\gamma} and \msosc{\gamma} for large OPF problems
representing Poland (PL)~\cite{matpower} and other European networks
from the \mbox{PEGASE} project~\cite{pegase} for which the first-order
relaxation fails to yield the globally optimal decision variables. The
hierarchies as presented in this paper have difficulty solving test
cases that minimize generation costs, due to the need for higher-order
constraints at many buses.\footnote{For problems that minimize generation cost, see~\cite{cdc2015} for a related approach that finds feasible points with objective values near the global optimum.} We therefore chose to minimize
  active power losses; the relaxation hierarchies globally solve these
  problems with higher-order constraints at only a few buses.

% Discuss termination after a maximum solver time of 2 hours or something, with these values added into the results

The low-impedance line preprocessing method in~\cite{cdc2015} with a threshold of $1\times 10^{-3}$~per~unit for the Polish systems and $3\times 10^{-3}$~per~unit for the PEGASE systems was used to improve numerical convergence. To further improve numerical convergence and computational speed, we did not enforce the sphere constraint for \msosc{\gamma}. Low-order relaxations from the complex hierarchy converged for the test cases considered here despite the lack of a convergence guarantee in the absence of the sphere constraint~\cite{complex_hierarchy}. Bounds on the lifted variables $y$ and $\hat{y}$ derived from the voltage magnitude limits~\eqref{eq:opfV} are also enforced to improve numeric convergence. 

The relaxations were implemented using MATLAB~2013a, \mbox{YALMIP 2015.06.26}~\cite{yalmip}, and \mbox{Mosek 7.1.0.28} and were solved using a computer with a quad-core 2.70~GHz processor and 16 GB of RAM. The results do not include the typically small formulation times. We used a tolerance $\epsilon = 1$~MVA.

\begin{figure}[t]
\centering
\includegraphics[totalheight=0.2925\textheight]{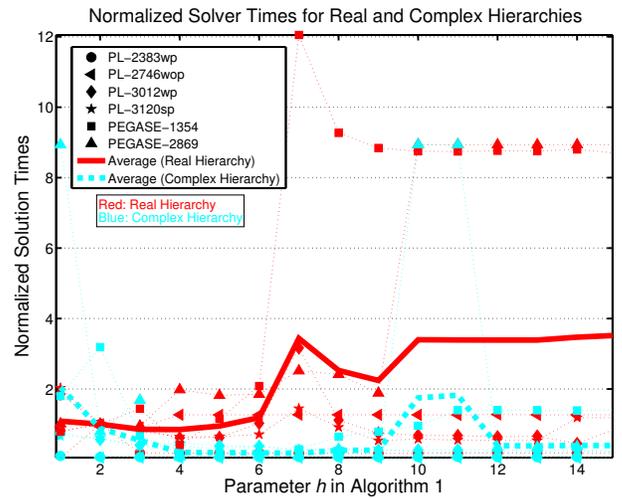}
\vspace{-5mm}
\caption{Solver times (normalized such that $h=2$ for the real
  hierarchy corresponds to a value of 1) vs. $h$} \label{f:tsolve}
\end{figure}

\begin{figure}[t]
\centering
\includegraphics[totalheight=0.2925\textheight]{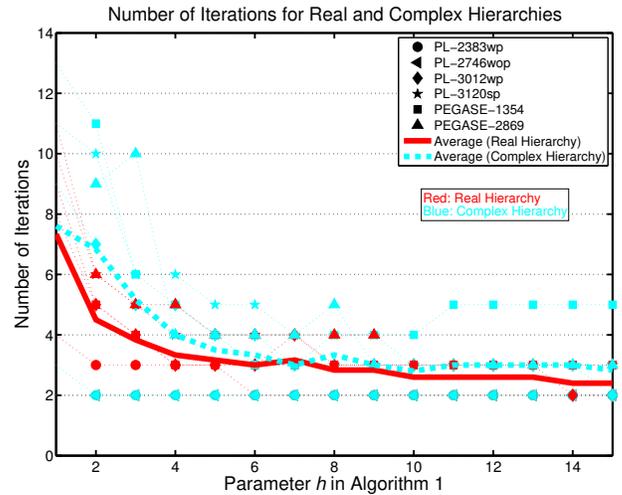}
\vspace{-5mm}
\caption{Number of Iterations vs. $h$} \label{f:niter}
\end{figure}

Figure~\ref{f:tsolve} shows a plot of the solver times for the large test cases vs. $h$. In order to match previous results~\cite{molzahn_hiskens-sparse_moment_opf,cdc2015,complex_hierarchy}, the values in Figure~\ref{f:tsolve} are normalized such that a value of one represents the solution time for the real hierarchy with $h = 2$.\footnote{A maximum solver time of two hours was enforced. Any convergence failures (e.g., insufficient memory) are assigned this maximum solver time. Solver times are reported as the mean time from three MOSEK solves.} The solver times for the real hierarchy with $h=2$ are 571, 2610, 261, 339, 351, and 806 seconds for the systems in the order shown in the legend of Figure~\ref{f:tsolve}. Figure~\ref{f:niter} shows a plot of the number of iterations of Algorithm~\ref{alg:heuristic} vs. $h$. 

%616, 2586, 267, 355, 345, and 805 seconds 

% update solver times once we have more results to average in.

%\footnote{See Tables~4 and 5 in~\cite{complex_hierarchy} for full results with $h=2$.}

%Table~\ref{t:h2_times} shows the results for $h=2$, including the number of iterations, globally optimal objective value, maximum mismatch at convergence $S_i^{\text{inj}}$, and solver times. All solver times can be derived from the results in Figure~\ref{f:tsolve} and Table~\ref{t:h2_times}.

These results demonstrate that the complex hierarchy is generally faster than the real hierarchy for a range of values of $h$ despite sometimes requiring more iterations. This is expected due to the smaller matrices in the complex hierarchy. Averaged over all test cases and choices of $h$, the complex hierarchy was a factor of 9.2 faster than the real hierarchy.
% 11.6 faster than the real hierarchy.

% Update factor of speed increase when I have all results averaged in

The computational challenges of \msosr{\gamma} are most apparent when Algorithm~\ref{alg:heuristic} increases $\gamma_i$ for buses that are contained within large cliques, which correspond to large matrices. For instance, \mbox{PEGASE-1354} with $h=7$ results in second-order relaxation constraints for a clique with ten buses and a solver time of 4112 seconds for this iteration of Algorithm~\ref{alg:heuristic}. This motivates future work in developing a more sophisticated algorithm that considers the size of the corresponding matrix rather than just the power injection mismatch when determining how to increment $\gamma_i$. This example also demonstrates the significant advantages of the complex hierarchy, for which global solution with the same set of higher-order buses only requires 89 seconds. With smaller matrices associated with the higher-order constraints, the complex hierarchy is also generally less sensitive to the parameter $h$ than the real hierarchy.

% as evidenced by the consistently fast solution times in Figure~\ref{f:tsolve}.

% \footnote{Ten-bus systems are themselves at the edge of tractability.}

% for \mbox{PEGASE-1354}, $h=7$ has 238   369   640, which heuristic 6 doesn't have. Heuristic 6 has 671, which heuristic 7 doesn't have. What is the story with these buses? Does one of them live in some large clique? Does the complex hierarchy not choose that bus or is it just better at handling the larger matrices?
% <-- Bus 640 is contained in a clique with 10 buses. That get bumped up to second-order and kills the solver times. Apply the same order vector to the complex hierarchy works much faster (89.3254 seconds, gives a solution within the tolerance 0.7165 MVA, vs. 4.112435e+03 seconds for that iteration of the real hierarchy, maxmismatch =0.818127 MVA)

Figure~\ref{f:niter} numerically demonstrates the expected trade-off related to $h$: with more higher-order constraints added each iteration, there are generally fewer iterations required but more time per iteration. The number of iterations quickly decreases but flattens out when $h \approx 4$. This flattening indicates a saturation in the benefit associated with adding more higher-order constraints at each iteration. For the loss minimization test cases, the average solver time results suggest that good values of $h$ are $7$ for \msosc{\gamma} and $4$ for \msosr{\gamma}.

% verify after averaging in more solver times

The results may be sensitive to the test cases and implementation details. Future work includes studying related hierarchies (e.g., alternatives to Algorithm~\ref{alg:heuristic}, the ``moment+penalization'' approach in~\cite{cdc2015}, and the mixed SDP/SOCP hierarchy in~\cite{powertech2015}) as well as other means of tightening the relaxations (e.g., the valid inequalities in~\cite{sun2015}). Leveraging efforts in test case development~\cite{griddata}, future work also includes extending this study using additional test cases.
% wang2015

\section{Conclusion}
\label{l:conclusion}
This paper presented two hierarchies of convex relaxations for globally solving OPF problems. The real hierarchy \msosr{\gamma} has been successfully applied to OPF problems in several previous works, while the complex hierarchy \msosc{\gamma} has only recently been introduced. After summarizing and comparing these hierarchies, we discussed an iterative method for exploiting sparsity and selectively applying computationally intensive higher-order relaxation constraints using a ``power injection mismatch'' approach to identify problematic buses. This approach requires a single parameter $h$ defining the maximum number of buses that are assigned higher-order constraints at each iteration of the algorithm. We presented numerical experiments using several large OPF problems to identify the impact of different choices for $h$.

% The 'triggered' command can be changed if desired:



\bibliographystyle{IEEEtran}
\bibliography{pscc2016}{}


%\begin{thebibliography}{1}
%\bibitem{Shell}
%M.~Shell, \emph{How to Use the IEEEtran Latex Class}, Latex Archive Contents, \verb+http://www.ieee.org/conferences_events/+ \verb+conferences/publishing/templates.htm+

%\end{thebibliography}

\end{document}